\documentclass[10pt,twocolumn,twoside]{IEEEtran} 

\IEEEoverridecommandlockouts                              % This command is only needed if 
\usepackage{amsmath,amsfonts,amssymb,mathptmx,bm}
\usepackage{diagbox,graphicx,subfigure,lineno,cases,cite}
\usepackage{epstopdf,pgf,tikz,color,xcolor}
\usepackage{graphicx,subfigure,pgf,tikz,pst-node}
\usetikzlibrary{arrows,shapes,automata,positioning,fit}
\allowdisplaybreaks % allow the equation display in different pages
\newcommand{\lm}{\lambda}
\newcommand{\e}{\varepsilon}
\newcommand{\R}{\mathbb{R}}

\newtheorem{thm}{Theorem}

\newtheorem{lem}{Lemma}
\newtheorem{ass}{Assumption}

\newtheorem{rem}{Remark}
\renewcommand{\bar }{\overline}
\newcommand{\pb}{\begin{IEEEproof} }
\newcommand{\pe}{\end{IEEEproof}}

\begin{document}
	\title{\LARGE \bf  Multi-agent Optimal Consensus with Unknown Control Directions \thanks{This work was supported by National Natural Science Foundation of China under Grant 61973043. }
	}
	\author{Yutao Tang \thanks{Y. Tang is with the School of Automation, Beijing University of Posts and Telecommunications, Beijing 100876, China (e-mail: yttang@bupt.edu.cn).}
	}
	\maketitle

	\begin{abstract}
		This paper studies an optimal consensus problem for a group of heterogeneous high-order agents with unknown control directions.   Compared with existing consensus results, the consensus point is further required to an optimal solution to some distributed optimization problem. To solve this problem, we first augment each agent with an optimal signal generator to reproduce the global optimal point of the given distributed optimization problem, and then complete the global optimal consensus design  by developing some adaptive tracking controllers for these augmented agents. Moreover, we present an extension when only real-time gradients are available. The trajectories of all agents in both cases are shown to be well-defined and achieve the expected consensus on the optimal point.  Two numerical examples are given to verify the efficacy of our algorithms.
	\end{abstract}
%	
%	\begin{IEEEkeywords}
%		Optimal consensus, high-order dynamics, unknown control direction, adaptive control
%	\end{IEEEkeywords}
\IEEEpeerreviewmaketitle

\section{Introduction}

Consensus is a fundamental problem in the field of multi-agent coordination and has been actively studied for decades. As an extended version of pure consensus, optimal consensus has been paid more and more attention due to its wide applications in multi-robot networks, machine learning, and big data technologies.  In the optimal consensus problem, each agent has a local cost function and all agents are expected to reach a consensus state that minimizes the sum of these individual cost functions. Many effective algorithms have been proposed for single-integrator multi-agent systems to achieve this optimal consensus goal under various conditions (see \cite{nedic2018distributed,yang2019survey} and references therein). 

Along with these fruitful optimal consensus results for single-integrator multi-agent systems, there are numerous optimal consensus tasks implemented by or depending on engineering multi-agent systems of high-order dynamics, e.g., source seeking in mobile sensor networks \cite{zhang2011extremum},  frequency control in power systems \cite{zhangxuan2017distributed}, and attitude formation control of rigid bodies \cite{song2017relative}. Thus, many authors seek to solve the optimal consensus problem for non-single-integrator multi-agent systems. Some recent attempts have been made for  second-order ones \cite{zhang2017distributed,xie2017global,qiu2019distributed}, general linear ones \cite{tang2019cyb}, and several classes of nonlinear multi-agent systems \cite{wang2016distributed,tang2020optimal}. 

So far, all these optimal consensus works  were only devoted to the cases when we have a prior knowledge of the control directions of agents' dynamics. Note that the control direction  of an engineering  plant may not always be known beforehand. Even it is known at first, it could be changed by some structural damages in many applications as shown in \cite{du2007adaptive, liu2008multivariable}. Therefore, it is crucial to consider the unknown control direction issue when resolving the optimal consensus problem for high-order engineering multi-agent systems.

At the same time, a plenty of pure consensus results without such optimization requirements have been derived for multi-agent systems with unknown control directions from integrators to nonlinear ones even with uncertainties by extending the classical Nussbaum-type controls \cite{nussbaum1983some,ye1998adaptive} to  decentralized and distributed cases, e.g., \cite{peng2014cooperative, chen2013adaptive, guo2016cooperative, rezaei2018adaptive, huang2018fully}. It is thus very interesting to ask whether similar Nussbaum-type controls can be constructed to tackle the optimal consensus problem in the presence of unknown control directions.

%For example,  \cite{peng2014cooperative} proposed a  Nussbaum-type adaptive controller for single-integrator agent such that consensus for this multi-agent system can be achieved. Parameterized uncertainties were also considered for both first and second order agents by constructing a special type of Nussbaum-type functions in \cite{chen2013adaptive}. However, this result relied on the assumption that all high-frequency gains have the same sign. This limitation has been removed in \cite{guo2016cooperative} for multi-agent systems with more general dynamics. Nevertheless, when the control directions are unknown, the solvability of optimal consensus problem even for single-integrator agents is still unclear.

Based on the aforementioned observations, we consider a group of high-order multi-agent systems with unknown control directions and seek distributed rules to solve the associated optimal consensus problem. Although some interesting optimal consensus results are available for these multi-agent systems \cite{shi2013reaching, zhang2017distributed, xie2017global, tang2019cyb}, the solvability of optimal consensus for them in the presence of unknown control directions is much more challenging and is still unclear. In fact, the gradient-based rules are basically nonlinear in light of the optimization requirement for the multi-agent system. More importantly, the unknown control directions of these agents and heterogeneous system orders  bring many extra technical difficulties to the associated optimal consensus analysis and design.

Motivated by the given designs in \cite{tang2019cyb}, we aim to develop an embedded control to solve the formulated optimal consensus problem for these agents. We will first assume the local function's analytic form and augment each agent with an optimal signal generator to reproduce the global optimal solution. Then, the expected optimal consensus will be achieved by embedding this generator into a Nussbaum-type adaptive tracking controller for each agent. Next, we will present an extension of the preceding designs using only real-time gradient information to achieve this optimal consensus.

The main contribution of this paper can be summarized as follows. First, compared with existing optimal consensus results assuming the knowledge of the control directions \cite{xie2017global,qiu2019distributed,tang2019cyb}, we remove this requirement and present effective distributed controllers for these heterogeneous high-order agents to reach an optimal consensus in the presence of unknown control directions.  To our knowledge, no other work solves such an optimal consensus problem under these circumstances yet. Second, as pure/average consensus can be achieved by solving some special optimal consensus problem, our algorithms naturally provide an alternative way other than \cite{chen2013adaptive,peng2014cooperative,rezaei2018adaptive,huang2018fully} to tackle such pure and average problems for agents with unknown control directions extending the derived consensus results in \cite{ren2008distributed,rezaee2015average}.% other than similar consensus results derived in \cite{chen2013adaptive, peng2014cooperative}. 

The rest of this paper is organized as follows. Some preliminaries are provided in Section \ref{sec:pre}. The problem formulation part is given in Section \ref{sec:form}. Main results are presented in Section \ref{sec:main}. Finally, simulations and our concluding remarks are presented at Sections \ref{sec:simu} and \ref{sec:con}.

\section{Preliminaries}\label{sec:pre}

We will use standard notations. Let $\R^N$ represent the $N$-dimensional Euclidean space.  Denote $||a||$ the Euclidean norm of a vector $a$ and $||A||$ the spectral norm of a matrix $A$.  ${\bm 1}_N$ (or ${\bm 0}_N$) denotes an $N$-dimensional all-one (or all-zero) column vector, and $I_{N}$ denotes the $N$-dimensional identity matrix.   Let $M_1=\frac{1}{\sqrt{N}} {\bm 1}_N$ and $M_2$ be the matrix satisfying $M_2^\top M_1={\bm 0}_{N-1}$, $M_2^\top M_2=I_{N-1}$ and $M_2 M_2^\top=I_{N}-M_1 M_1^\top$.  We may omit the subscript when it is self-evident.

A directed graph (digraph) is described by $\mathcal{G}=(\mathcal{N}, \,\mathcal{E}, \, \mathcal{A})$ with node set $\mathcal{N}=\{1,{\dots},N\}$ and edge set $\mathcal {E}$. $(i,\,j)\in \mathcal{E}$ denotes an edge from node $i$ to $j$. The weighted adjacency matrix $\mathcal{A}=[a_{ij}]\in \mathbb{R}^{N\times N}$ is defined by $a_{ii}=0$ and $a_{ij}\geq 0$. Here $a_{ij}>0$ iff $(j,\,i)\in \mathcal{E}$.  Node $i$'s neighbor set is defined as $\mathcal{N}_i=\{j\mid (j,\, i)\in \mathcal{E} \}$. %We denote $\mathcal{N}_i^0=\mathcal{N}_i\cup \{i\}$.  
A directed path is an ordered sequence of vertices such that each intermediate pair of vertices is an edge. %If $a_{ij}=a_{ji}$ for any $i,\,j\in  \mathcal{N}$, we say this graph is undirected. 
If there is a directed path between any two nodes, then the digraph is said to be strongly connected. %A strongly connected undirected graph is said to be connected.  
The in-degree and out-degree of node $i$ are defined as $d^{\mbox{in}}_i=\sum\nolimits_{j=1}^N a_{ij}$ and $d^{\mbox{out}}_i=\sum\nolimits_{j=1}^N a_{ji}$. The Laplacian of digraph $\mathcal{G}$ is defined as $L\triangleq D^{\mbox{in}}-\mathcal{A}$ with $D^{\mbox{in}}=\mbox{diag}(d^{\mbox{in}}_1,\,\dots,\,d^{\mbox{in}}_N)$. A digraph is weight-balanced if $d^{\mbox{in}}_i=d^{\mbox{out}}_i$ for any $i\in \mathcal{N}$.  Note that $L{\bm 1}_N={\bm 0}_N$ for any digraph. %If the digraph is strongly connected, $0$ is a simple eigenvalue of its Laplacian.  
A digraph is weight-balanced iff ${\bm 1}_N^\top L={\bm 0}_N$, which is also equivalent to ${\mbox{Sym}(L)=\frac{L+L^\top }{2}}$ being positive semidefinite.  For a strongly connected and weight-balanced digraph, %$0$ is a simple eigenvalue of $\mbox{Sym}(L)$ and all other eigenvalues are positive real numbers. In this case, 
we can order the eigenvalues of $\mbox{Sym}(L)$ as $0=\lambda_1<\lambda_2\leq \dots\leq \lambda_N$ and  have $ \lambda_2 I_{N-1}\leq M_2^\top \mbox{Sym}(L)M_2\leq \lambda_N I_{N-1}$. See \cite{godsil2001algebraic} for more details. 

%\subsection{Convex analysis}\
A function $f\colon \R^m \rightarrow \R $ is said to be convex if $f(a\zeta_1+(1-a)\zeta_2)\leq af(\zeta_1)+(1-a)f(\zeta_2)$ for any  $0\leq a \leq 1$ and all $\zeta_1,\,\zeta_2 \in \R^m$. When $f$ is differentiable,  it is convex if $f(\zeta_1)-f(\zeta_2)\geq \nabla f(\zeta_2)^\top (\zeta_1 -\zeta_2)$ for all $\zeta_1,\,\zeta_2 \in \mathbb{R}^m$. We say $f$ is $\omega$-strongly convex  over $\R^m$ if $[\nabla f(\zeta_1)-\nabla f(\zeta_2)]^\top (\zeta_1 -\zeta_2)\geq \omega ||\zeta_1 -\zeta_2||^2$ for all  $\zeta_1,\,\zeta_2 \in \R^m$ with $\omega >0$. A vector-valued function ${\bm f}\colon \R^m \rightarrow \R^m$ is said to be $\vartheta$-Lipschitz  if $ ||{\bm f}(\zeta_1)-{\bm f}(\zeta_2)||\leq \vartheta ||\zeta_1-\zeta_2||$ for all  $\zeta_1, \,\zeta_2 \in \R^m$ with $\vartheta>0$.

\section{Problem Formulation}\label{sec:form}
Consider a heterogeneous multi-agent system consisting of $N$ agents described by
\begin{align}\label{sys:agent}
{y}^{(n_i)}_i=b_i u_i,\quad i\in \mathcal{N}\triangleq \{1,\,\dots,\,N\}
\end{align}
where $y_i\in \mathbb{R}$ and $u_i\in \R$ are its output and input, respectively. Integer $n_i\geq 1$ is the order of system \eqref{sys:agent} and constant $b_i$ is assumed to be away from zero but unknown. This constant $b_i$ is often called the high-frequency gain of agent \eqref{sys:agent}, which represents the motion direction of this agent in any control strategy. The parameters $n_i$ and $b_i$ of each agent are allowed to be different from each other.  

We endow agent $i$ with a local cost function $f_i\colon\R\to\R$ for $i\in \mathcal{N}$, and define the global cost function as $f(y)=\sum_{i=1}^{N} f_i(y)$. For multi-agent system \eqref{sys:agent}, we aim to develop an algorithm such that all agent outputs achieve a consensus on the minimizer to this global cost function. 

The following assumption is often made in  literature \cite{ruszczynski2011nonlinear,jakovetic2015linear,kia2015distributed,tang2020optimal}, which guarantees the existence and uniqueness of the minimal solution to function $f$. 
\begin{ass}\label{ass:convexity-strong}
	For $i\in \mathcal{N}$,  function $f_i $ is $\underline{l}$-strongly convex and its gradient is $\overline l$-Lipschitz for two constants $\bar{l}\geq \underline l>0$.
\end{ass}

As usual, we assume this optimal solution is finite and denote it as $y^{\star}$, i.e.
\begin{align}\label{opt:main}
y^{\star}={\arg}\min_{ y \in \R} \; f(y)=\sum\nolimits_{i=1}^{N} f_i(y)
\end{align}

Due to the privacy of local cost function $f_i$,  no agent can unilaterally  determine the global optimal solution $y^{\star}$ by itself. Hence, our problem cannot be solved without cooperation and information sharing among these agents.  For this purpose, we use a digraph $\mathcal {G}=(\mathcal {N}, \mathcal {E}, \mathcal{A})$ to describe the information sharing topology. An edge $(j,\,i )\in \mathcal{E}$ with weight $a_{ij}>0$ means that agent $i$ can get the information of agent $j$.

To guarantee that any agent's information can reach any other agents,  we suppose the following assumption holds. % as that in many publications \cite{tang2015distributed,yi2015distributed}.
\begin{ass}\label{ass:graph}
	The digraph $\mathcal{G}$ is weight-balanced  and strongly connected.
\end{ass}
%This assumption is about the connectivity of information sharing graph $\mathcal{G}$, which guarantees that any agent's information can reach any other agents. 

Then, our optimal consensus problem is to design $u_i$ for agent $i$ under the information constraint described by digraph $\mathcal{G}$, such that these agents achieve an optimal consensus determined by the global objective function $f$ in the sense that $y_i-y^{\star}\to 0$ as $t \to \infty$ for any $i\in \mathcal{N}$, while the trajectories of this multi-agent system are maintained to be bounded.

\begin{rem}\label{rem:form}
	This optimal consensus problem has been extensively studied in literature for multi-agent systems assuming the high-frequency gain is known \cite{zhang2017distributed,xie2017global,qiu2019distributed,tang2019cyb}. But in our work, this prior knowledge of each agent’s control direction is no longer necessary, which means that agents may have different and unknown control directions. To the best of our knowledge, no other works have studied the optimal consensus problem under these circumstances yet.
\end{rem}

It is interesting to point out that when the local cost functions are chosen as $f_i(y)=c_i(y-y_i(0))^2$ with $c_i>0$ for each $i \in \mathcal{N}$, we can solve a scaled consensus problem with the final consensus point $y^{\star}=\frac{\sum_{i=1}^N c_i y_i(0)}{\sum_{i=1}^N c_i}$.  Thus, this formulation  provides an applicable way to solve their pure and average consensus problems  for these high-order agents in the presence of unknown control directions.

\section{Main Result}\label{sec:main}

In this section, we will present an embedded design  to solve our formulated optimal consensus problem following the technical line developed in \cite{tang2019cyb}.

To this end, we first consider an auxiliary optimal consensus problem with the same requirement for agents in form of $\dot{r}_i=\mu_i$ and then convert our problem into an output tracking control problem for agent \eqref{sys:agent} with reference $r_i$. As the former subproblem is essentially a conventional optimal consensus problem for single-integrator multi-agent system with $b_i=1$ and has been well-studied in existing literature, we use the following optimal signal generator to complete our design:
\begin{align}\label{sys:generator}
\dot{r}_i&=-\alpha \nabla f_i(r_i)-\beta \sum\nolimits_{j=1}^{N}a_{ij}(r_i-r_j)-\sum\nolimits_{j=1}^{N}a_{ij}(v_i-v_j)\nonumber\\
\dot{v}_i&=\alpha \beta  \sum\nolimits_{j=1}^{N}a_{ij}(r_i-r_j)
\end{align}
where $\alpha,\,\beta$ are constants to be specified later.  Putting it into a compact form gives
\begin{align}\label{sys:composite-osg}
\begin{split}
\dot{r}&=-\alpha \nabla \tilde f(r)- \beta Lr-Lv,\quad \dot{v}=\alpha \beta Lr
\end{split}
\end{align}
where $r=\mbox{col}(r_1,\,\dots,\,r_N)$, $v=\mbox{col}(v_1,\,\dots,\,v_N)$, and $\tilde f(r)\triangleq \sum\nolimits_{i=1}^Nf_i(r_i)$ is $\underline{l}$-strongly convex while its gradient $ \nabla\tilde f(r)$ is $\bar l$-Lipschitz with $\bar l=\max_i\{\bar l_i\}$ and $\underline{l}=\min_i\{\underline{l}_i\}$. 

System \eqref{sys:generator} is a distributed primal-dual variant to determine the optimal consensus point $y^{\star}$. Its effectiveness has already been established in \cite{tang2020optimal} by semistability arguments. Here, we provide a sketch of proof using Lyapunov stability analysis.
\begin{lem}\label{lem:generator}
		Suppose Assumptions \ref{ass:convexity-strong}--\ref{ass:graph} hold and let $\alpha\geq \max\{1,\,\frac{1}{\underline{l}},\,\frac{2\bar l^2}{\underline{l}\lm_2}\}$, $\beta\geq \max\{1,\, \frac{1}{\lambda_2},\,\frac{6\alpha^2\lambda_N^2}{\lambda_2^2} \}$.  Then,  the trajectory of system \eqref{sys:generator} from any initial point is bounded over $[0,\,\infty)$ and $r_i(t)$ approaches  $y^{\star}$ exponentially as $t\to \infty$ for $i\in\mathcal{N}$. 
\end{lem}
\pb The basic idea is to perform a change of coordinates and determine a reduced-order system with a unique equilibrium point to avoid semistability arguments. 

Let $\mbox{col}(r^{\star},\,v^{\star})$ be any equilibrium point of system \eqref{sys:composite-osg} and can verify $r^{\star}={\bm 1}_N y^{\star}$ under Assumptions \ref{ass:convexity-strong}--\ref{ass:graph} by Theorem 3.27 in \cite{ruszczynski2011nonlinear}.  Perform the coordinate transformation: $\bar r_1=M_1^\top (r-r^\star)$, $\bar r_2=M_2^\top (r-r^\star)$, $\bar v_1=M_1^\top (v-v^\star)$, and $\bar v_2=M_2^\top[( v+\alpha r)-( v^\star+\alpha r^\star)]$. It follows that $\dot{\bar v}_1=0$ and  
\begin{align}\label{sys:composite-osg-reduced}
\begin{split}
\dot{\bar r}_1&=-\alpha M_1^\top {\bm \Pi}\\
\dot{\bar r}_2&=-\alpha M_2^\top {\bm \Pi}-\beta M_L \bar r_2 + \alpha M_L \bar r_2-M_L\bar v_2\\
\dot{\bar v}_2&=-\alpha M_L {\bar v}_2+\alpha^2 M_L \bar r_2-\alpha^2 M_2^\top {\bm \Pi} 
\end{split}
\end{align}
where  ${\bm \Pi}\triangleq \nabla \tilde f(r)-\nabla \tilde f(r^\star)$ and $M_L= M_2^\top LM_2$.   Let $\bar r=\mbox{col}(\bar r_1,\,\bar r_2)$, and $V_{\rm o}(r,\, v)=\bar r^\top \bar r + \frac{1}{\alpha^3} \bar v_1^\top \bar v_1+\frac{1}{\alpha^3}\bar v_2^\top \bar v_2$ in this new coordinate. It is quadratic and positive definite. By Young's inequality to handle the cross terms as that in \cite{tang2020optimal}, the derivative of $V_{\rm o}(t)$ along the trajectory of \eqref{sys:composite-osg} satisfies 
\begin{align*}
\dot{V}_{\rm o}%&=-2\alpha (r-r^\star)^\top {\bm \Pi}+ 2\bar r_2^\top [-\beta M_L \bar r_2 + \alpha M_L \bar r_2-M_L\bar v_2 ]\\
%&+\frac{2}{\alpha^3}\bar v_2^\top[-\alpha M_L {\bar v}_2+\alpha^2 M_L \bar r_2-\alpha^2 M_2^\top {\bm \Pi} ]\\
%\leq & -2\alpha\underline{l}||\bar r||^2-2\beta\lambda_2||\bar r_2||^2+2\alpha  \lambda_N  ||\bar r_2||^2+2 \lambda_N  ||\bar r_2||||\bar v_2||\\
%&-\frac{2\lambda_2}{\alpha^2}||\bar v_2||^2+ \frac{2}{\alpha} \lambda_N ||\bar r_2|| ||\bar v_2||+\frac{2\bar l}{\alpha}||\bar v_2||||\bar r||\\
%%\leq\;& -2\alpha\underline l ||\bar r||^2-2\beta\lambda_2||\bar r_2||^2+2\alpha  \lambda_N  ||\bar r_2||^2+\frac{\lambda_2}{3\alpha^2}||\bar v_2||^2+\frac{3\alpha^2\lambda_N^2}{\lambda_2}||\bar r_2||^2\\
%&-\frac{2\lambda_2}{\alpha^2}||\bar v_2||^2+ \frac{\lambda_2}{3\alpha^2}||\bar v_2||^2+\frac{3\lambda_N^2}{\lambda_2}||\bar r_2||^2+\frac{\lambda_2}{3\alpha^2}||\bar v_2||+\frac{3\bar l^2}{\lambda_2}||\bar r||^2\\
%\leq& -(2\alpha\underline{l} -\frac{3\bar l^2}{\lambda_2} )||\bar r||^2-\frac{\lambda_2}{\alpha^2}||\bar v_2||^2\\
%& - (2\beta\lambda_2-2\alpha\lambda_N-\frac{3\alpha^2\lambda_N^2}{\lambda_2}-\frac{3\lambda_N^2}{\lambda_2})||\bar r_2||^2\\
&\leq - \frac{1}{2}||\bar r||^2-\frac{1}{2\alpha^3}||\bar v_2||^2\triangleq W_{\rm o}(\bar r,\, \bar v_2)
\end{align*}
This implies the Lyapunov stability of system \eqref{sys:composite-osg} at $\mbox{col}(r^{\star},\,v^{\star})$ and the boundedness of the trajectory from any initial point over $[0,\,\infty)$. Further considering the reduced-order system \eqref{sys:composite-osg-reduced} with a Lyapunov function $W_{\rm o}(\bar r,\,\bar v_2)$,  one can obtain that $\dot{W}_{\rm o}\leq -\frac{1}{2}W_{\rm o}$ along the trajectory  of \eqref{sys:composite-osg-reduced}. Recalling Theorem 4.10 in \cite{khalil2002nonlinear},  $W_{\rm o}(\bar r(t),\,\bar v_2(t))$ and  $\bar r(t)$ exponentially converge to $0$ as $t$ goes to infinity. The proof is complete.
\pe

With this generator  \eqref{sys:generator}, each agent can get an asymptotic estimate $r_i$ of the global optimizer $y^{\star}$. Thus, we are left to solve an output tracking problem for agent $i$ with reference $r_i$. 

When $b_i=1$, a pole-placement based tracking control was presented in \cite{tang2019cyb} for  multi-agent system \eqref{sys:agent} to complete the whole design. Controllers with bounded constraints were also developed to achieve an optimal consensus in literature \cite{xie2017global,qiu2019distributed}. However, the control directions are assumed to be unknown in our current case. Consequently,  such methods are no longer applicable to agent \eqref{sys:agent} and we have to seek new tracking rules to solve our optimal consensus problem.
 
For this purpose, we denote $y_{i1}=y_i-r_i$ and $y_{i\iota}\triangleq \epsilon^{\iota-1} y_{i}^{(\iota-1)}$ for $2\leq \iota\leq n_i$ with a constant  $\epsilon>0$ to be specified later. Choose constants $k_{i\iota}$ for $1\leq \iota\leq n_i-1$ such that the polynomial $p_i(\lambda)=\sum_{\iota=1}^{n_i-1}k_{i\iota} \lambda^{\iota-1}+\lambda^{n_i-1}$ is Hurwitz. Letting $z_i=\mbox{col}(y_{i1},\,\dots,\, y_{in_i-1})$ and $\zeta_i=\sum_{\iota=1}^{n_i-1}k_{i\iota} y_{i\iota}+y_{in_i}$  gives the following translated multi-agent system: 
\begin{align}\label{sys:rd=1}
\begin{split}
\dot{z}_i&=\frac{1}{\epsilon}A_{i1}z_i+ \frac{1}{\epsilon} A_{i2}\zeta_i+E_{i1}\dot{r}_i\\
\dot{\zeta}_i&=\frac{1}{\epsilon}A_{i3} z_i+\frac{1}{\epsilon}A_{i4} \zeta_i+ \epsilon^{n_i-1} b_i u_i+ E_{i2}\dot{r}_i
\end{split}
\end{align}
where the associated matrices are defined as follows.
\begin{align*}
&A_{i1}=\left[\begin{array}{c|c}
{\bm 0}_{n_i-2}&I_{n_i-2}\\\hline
-k_{i1}&-k_{i2}~~\cdots~~-k_{i~n_i-1}
\end{array}\right], \quad A_{i2}=\begin{bmatrix}
{\bm 0}_{n_i-2}\\
1
\end{bmatrix}\\
&A_{i3}=\begin{bmatrix} -k_{in_i-1}k_{i1}~~k_{i1}-k_{in_i-1}k_{i2}~~\cdots~~k_{i n_i-2}-k_{in_i-1}^2
\end{bmatrix}\\
&A_{i4}=k_{in_i-1},\quad E_{i1}=\begin{bmatrix}
1~~{\bm 0}_{n_i-2}^\top 
\end{bmatrix}^\top ,\quad E_{i2}=-k_{i1}
\end{align*}

%Note that the righthand size of \eqref{sys:composite-osg} is global Lipschitz under Assumption \ref{ass:convexity-strong}. Jointly with the exponetially, $\dot{r}_i$ exponentially vanishes as $t$ goes to $\infty$.  
From the proof of Lemma \ref{lem:generator},  system \eqref{sys:composite-osg-reduced} is exponentially stable at the origin. Joitly with the Lipschitzness of $\Pi$ under Assumption \ref{ass:convexity-strong} and $\dot{r}=M_1\dot{\bar r}_1+M_2\dot{\bar r}_2$,  $\dot{r}(t)$ will exponentially converge to ${\bm 0}$. Thus, we have converted the original optimal consensus problem into a robust stabilization problem of the translated system \eqref{sys:rd=1} with time-decaying disturbances $\dot{r}_i$.

Motivated by the designs in \cite{ye1998adaptive,peng2014cooperative,chen2019nussbaum}, we use the following Nussbaum-type rule to serve the tracking purpose:
\begin{align*}
u_i&=\bar{\mathcal{N}}(\theta_i)\zeta_i,\quad \dot{\theta}_i=\zeta^2_i
\end{align*}
where % $\zeta_i$ is defined as above and 
$\bar{\mathcal{N}}$ is a smooth function satisfying:
\begin{align}\label{eq:nussbaum-type-A}
\limsup_{\theta \to\infty}\frac{\int_{0}^\theta \bar{\mathcal{N}}(s) {\rm d}s}{\theta}=\infty,\quad \liminf_{\theta \to\infty}\frac{\int_{0}^\theta \bar{\mathcal{N}}(s) {\rm d}s}{\theta}=-\infty
\end{align}
Commonly used examples include  $\theta^2\sin \theta$ and $e^{\theta^2}\sin \theta$.

The overall controller to solve our problem is then:
\begin{align}\label{ctr:main}
\begin{split}
u_i&=\bar{\mathcal{N}}(\theta_i)\zeta_i\\
\dot{\theta}_i&=\zeta^2_i\\
\dot{r}_i&=-\alpha \nabla f_i(r_i)-\beta \sum\nolimits_{j=1}^{N}a_{ij}(r_i-r_j)-\sum\nolimits_{j=1}^{N}a_{ij}(v_i-v_j)\\
\dot{v}_i&=\alpha \beta  \sum\nolimits_{j=1}^{N}a_{ij}(r_i-r_j)
\end{split}
\end{align}
 where $\zeta_i=k_{i1}(y_i-r_i)+\sum_{\iota=2}^{n_i-1}k_{i\iota} \epsilon^{\iota-1} y_{i}^{(\iota-1)} + \epsilon^{n_i-1} y_{i}^{(n_i-1)}$ defined as above. This controller is indeed distributed in the sense of using only agent $i$'s own and neighboring information.

It is time to present our first main theorem of this paper. 
\begin{thm}\label{thm:main}
	Consider the multi-agent system consisting of $N$ agents given by \eqref{sys:agent}. Suppose Assumptions \ref{ass:convexity-strong}--\ref{ass:graph} hold. Then, there exist two positive constants $ \alpha,\, \beta$ such that the optimal consensus problem for this multi-agent system \eqref{sys:agent} and \eqref{opt:main} is solved by the controller \eqref{ctr:main} for any $\epsilon>0$.
\end{thm} 
\pb According to Lemma \ref{lem:generator}, it suffices for us to solve the tracking problem for each agent, which can be further converted to the robust stabilization problem for the translated agent \eqref{sys:rd=1}. Hence, we only have to show the trajectory of the translated system \eqref{sys:agent} from any initial point is well-defined over the time interval $[0,\,\infty)$ and $y_{i1}$ converges to zero.

To this end, we first show that the trajectory of this multi-agent system is well-defined over the time interval $[0,\,\infty)$. Note that the local error system for agent $i$ is 
\begin{align*}
\dot{z}_i&=\frac{1}{\epsilon}A_{i1}z_i+ \frac{1}{\epsilon}A_{i2}\zeta_i+E_{i1}\dot{r}_i\\
\dot{\zeta}_i&=\frac{1}{\epsilon}A_{i3} z_i+\frac{1}{\epsilon}A_{i4} \zeta_i+\epsilon^{n_i-1} b_i\bar{\mathcal{N}}(\theta_i)\zeta_i+E_{i2}\dot{r}_i\\
\dot{\theta}_i&=\zeta^2_i\\
\dot{r}_i&=-\alpha \nabla f_i(r_i)-\beta \sum\nolimits_{j=1}^{N}a_{ij}(r_i-r_j)-\sum\nolimits_{j=1}^{N}a_{ij}(v_i-v_j)\nonumber\\
\dot{v}_i&=\alpha \beta  \sum\nolimits_{j=1}^{N}a_{ij}(r_i-r_j)
\end{align*}
where $A_{i1}$ is Hurwitz according to the choice of $k_{i\iota}$.  Thus, there must be a positive definite matrix $P_i\in \R^{n_i-1\times n_i-1}$ such that $A_{i1}^\top P_{i}+P_iA_{i1}=-2I_{n_i-1}$. From the smoothness of related functions, the trajectory of each subsystem is well-defined on its maximal interval $[0,\,t_{if})$.  We claim that $t_{if}=\infty$ for each $i$.  In the following, we will prove this by seeking a contradiction. 

Assume $t_{if}$ is finite. We are going to prove that all involved signals  are bounded over the time interval $[0,\,t_{if})$. Take $V_i(z_i,\,\zeta_i)=z_i^\top P_i z_i+\frac{1}{2}\zeta_i^2$ as a sub-Lyapunov function for agent $i$. It is positive definite with a time derivative along the trajectory of the above error system as follows.
\begin{align}\label{eq:thm1:eq1}
\dot{V}_i&= 2 z_i^\top P_i [\frac{1}{\epsilon}A_{i1}z_i+ \frac{1}{\epsilon}A_{i2}\zeta_i+E_{i1}\dot{r}_i]\nonumber \\ 
&\quad +\zeta_i(\frac{1}{\epsilon}A_{i3} z_i+\frac{1}{\epsilon}A_{i4} \zeta_i+\epsilon^{n_i-1}b_i\bar{\mathcal{N}}(\theta_i)\zeta_i+E_{i2}\dot{r}_i) \nonumber\\ 
&\leq -\frac{2}{\epsilon}||z_i||^2+ \frac{1}{3\e}||z_i||^2+\frac{3}{\e}||P_i A_{i2}||^2\zeta_i^2+\frac{1}{3\e}||z_i||^2 \nonumber\\ 
&\quad +{3\e}||P_iE_{i1}||^2\dot{r}_i^2+\frac{1}{3\e}||z_i||^2+\frac{3}{\e}||A_{i3}||^2\zeta_i^2+\frac{1}{\e} |A_{i4}| \zeta_i^2 \nonumber\\ 
&\quad + \epsilon^{n_i-1} b_i\bar{\mathcal{N}}(\theta_i)\zeta_i^2+\frac{1}{\e}\zeta_i^2+{\e}||E_{i2}||^2\dot{r}_i^2 \nonumber \\ 
&=-\frac{1}{\e}||z_i||^2+( \epsilon^{n_i-1} b_i\bar{\mathcal{N}}(\theta_i)+C_{i \theta_1})\zeta_i^2+ C_{i\theta_2}\dot{r}_i^2 %\nonumber\\
%&=-\frac{1}{\e}||z_i||^2-( \epsilon^{n_i-1}  b_i\bar{\mathcal{N}}(\theta_i)-C_{i \theta_1})\dot{\theta}_i+ C_{i\theta_2}\dot{r}_i^2  %&\leq -(b_i\bar{\mathcal{N}}(\theta_i)-C_{i\theta_1})\dot{\theta}_i+ C_{i \theta_2}\dot{r}_i^2
\end{align}
where  we use Young's inequality to handle the cross terms with constants $C_{i\theta_1}=\frac{1}{\e}(3||P_i A_{i2}||^2+3||A_{i3}||^2+|A_{i4}|+1)$ and $C_{i \theta_2}=3\e||P_i E_{i1}||^2+\e||E_{i2}||^2$. 

Recalling Lemma \ref{lem:generator},  $r_i(t)$ and $\dot{r}_i(t)$ exponentially converge to $y^{\star}$ and $0$ under Assumptions \ref{ass:convexity-strong}--\ref{ass:graph}. Thus, $\dot{r}_i(t)$ is square-integrable over $[0,\,\infty)$. Denote $V_i(t)\triangleq V_i(z_i(t),\,\zeta_i(t))$ for short. Noting that $\dot{\theta}_i=\zeta_i^2$, we integrate both sides of \eqref{eq:thm1:eq1} from $0$ to $t$ and have the following inequality for some constant $C_{i0}>0$:
\begin{align*}%\label{eq:thm1:eq2}
{V}_i(t)-V_i(0) \leq \epsilon^{n_i-1}b_i  \int_{\theta_i(0)}^{\theta_i(t)} \bar{\mathcal{N}}(s)\, {\rm d}s+C_{i \theta_1}\theta_i(t)+C_{i0}
\end{align*}
As $\theta_i(t)$ is monotonically increasing, it either has a finite limit or grows to $\infty$. Assuming $\theta_i(t)$ tends to $\infty$, we divide both sides by $\theta_i(t)$ for a large enough $t$ and have
\begin{align*}
0\leq \epsilon^{n_i-1}b_i \frac{ \int_{\theta_i(0)}^{\theta_i(t)} \bar{\mathcal{N}}(s)\, {\rm d}s}{\theta_i(t)}+C_{i \theta_1}+\frac{C_{i0}+V_i(0)}{\theta_i(t)} 
\end{align*}
According to the property \eqref{eq:nussbaum-type-A} of $\bar{\mathcal{N}}$, this inequality will finally be violated for any fixed $b_i$. Hence, $\theta_i(t)$ must be bounded over $[0, \,t_{if})$.  Recalling the controller \eqref{ctr:main}, $z_i(t)$, $\zeta_i(t)$,  $u_i(t)$, $\dot{\zeta}_i(t)$, and $\dot{\theta}_i(t)$ are also bounded over $[0, \,t_{if})$ for each $i\in \mathcal{N}$.  This implies with a contradiction that no finite-time escape phenomenon happens. Thus, we have $t_{if}=\infty$.

From the boundedness of $\dot{\theta}_i$, the function $\theta_i(t)$ is uniformly continuous with respect to time $t$. Note that 
\begin{align*}
\int_{0}^t \zeta_i^2(s){\rm d}s=\int_{0}^t  \dot{\theta}_i(s){\rm d} s \leq \theta_i(\infty)-\theta_i(0)
\end{align*}
Since $\theta_i(\infty)$ exists and is finite, $\zeta^2_i(t)$ is thus integrable. By Lemma 8.2 in \cite{khalil2002nonlinear}, we have $\zeta_i(t)\to 0$ as $t$ goes to $\infty$.  

Considering the $z_i$-subsystem, it is input-state stable with input $\frac{1}{\e}A_{i2}\zeta_i+E_{i1}\dot{r}_i$ and state $z_i$. Since both $\zeta_i(t)$ and $\dot{r}_i(t)$ converge to $0$ when $t$ goes to $\infty$, we recall Theorem 1 in \cite{sontag2003remark} and obtain that $|y_i(t)-r_i(t)|\to 0$ as $t$ goes to $\infty$. Jointly using the triangle inequality and the convergence of $r_i(t)$ to $y^{\star}$, we have that $|y_i(t)-y^{\star}|\leq|y_i(t)-r_i(t)|+|r_i(t)-y^{\star}|\to 0$ as $t$ goes to $\infty$. The proof is thus complete.
\pe

In controller \eqref{ctr:main}, we require  the analytic form of $\nabla f_i$ to ensure the feasibility of our optimal signal generator \eqref{sys:generator}. However, in many cases, only real-time gradient $\nabla f_i(y_i)$ is available for agent $i$ and the controller \eqref{ctr:main} is thus not implementable.

To tackle this issue, we limit us to the case when all high-frequency gains have the same sign.  Replacing $\nabla f_i(r_i)$ with the real-time gradient $\nabla f_i(y_i)$, we present the following controller:
\begin{align}\label{ctr:main-online}
u_i&=\bar{\mathcal{N}}(\theta_i)\zeta_i  \nonumber\\
\dot{\theta}_i&=\zeta^2_i\\
\dot{r}_i&=-\alpha \nabla f_i(y_i)-\beta \sum\nolimits_{j=1}^{N}a_{ij}(r_i-r_j)-\sum\nolimits_{j=1}^{N}a_{ij}(v_i-v_j) \nonumber\\
\dot{v}_i&=\alpha \beta  \sum\nolimits_{j=1}^{N}a_{ij}(r_i-r_j) \nonumber
\end{align}
where $\zeta_i$ is defined as in \eqref{ctr:main} and $\bar{\mathcal{N}}$ is strengthened to satisfy
\begin{align}\label{eq:pe-nussbaum-type-b}
\begin{split}
\lim_{\theta\to \infty}\frac{\int_{0}^{y}\bar{\mathcal{N}}_+(s){\rm d}\,s}{\theta}=\infty,\quad  \limsup_{\theta\to \infty}\frac{\int_{0}^{y}\bar{\mathcal{N}}_+(s){\rm d}\,s}{\int_{0}^{y}\bar{\mathcal{N}}_-(s){\rm d}\,s}=\infty\\
\lim_{\theta\to \infty}\frac{\int_{0}^{y}\bar{\mathcal{N}}_-(s){\rm d}\,s}{\theta}=\infty,\quad  \limsup_{\theta\to \infty}\frac{\int_{0}^{y}\bar{\mathcal{N}}_-(s){\rm d}\,s}{\int_{0}^{y}\bar{\mathcal{N}}_+(s){\rm d}\,s}=\infty
\end{split}
\end{align}
with $\bar{\mathcal{N}}_+(\theta)\triangleq\max\{\bar{\mathcal{N}}(\theta),\,0\}$, $\bar {\mathcal{N}}_-(\theta)\triangleq \max\{-\bar {\mathcal{N}}(\theta),\,0\}$ for any $\theta \in \R$.  
It can be verified that such functions satisfy the condition \eqref{eq:nussbaum-type-A} and thus are special Nussbaum functions. Some feasible examples have been used in literature \cite{ding2015adaptive,chen2019nussbaum}. 

\begin{thm}\label{thm:main-online}
	Consider the multi-agent system consisting of $N$ agents given by \eqref{sys:agent}. Suppose all high-frequency gains are unknown but with the same sign and Assumptions \ref{ass:convexity-strong}--\ref{ass:graph} hold. Then, there exist constants $\alpha,\,\beta>0$ and $\epsilon^{\star}>0$ such that the optimal consensus problem for this multi-agent system \eqref{sys:agent} and \eqref{opt:main} is solved by the controller \eqref{ctr:main-online} for any $0< \epsilon < \epsilon^{\star}$.
\end{thm}
\pb Basically, we will decrease the parameter $\e$ to compensate the difference between $\nabla f_i(r_i)$ and $\nabla f_i(y_i)$. 

By the proof of Lemma \ref{lem:generator}, The composite system in this case can be written as follows:
\begin{align*}%\label{sys:thm:main-online}
\dot{z}_i&=\frac{1}{\epsilon}A_{i1}z_i+ \frac{1}{\epsilon}A_{i2}\zeta_i+E_{i1}\dot{r}_i\\
\dot{\zeta}_i&=\frac{1}{\epsilon}A_{i3} z_i+\frac{1}{\epsilon}A_{i4} \zeta_i+\epsilon^{n_i-1} b_i\bar{\mathcal{N}}(\theta_i)\zeta_i+E_{i2}\dot{r}_i\\
\dot{\theta}_i&=\zeta^2_i\\
\dot{\bar r}_1&=-\alpha M_1^\top ({\bm \Pi}+{\bm \Pi}_1)\\
\dot{\bar r}_2&=-\alpha M_2^\top( {\bm \Pi}+{\bm \Pi}_1)-\beta M_L \bar r_2 + \alpha M_L \bar r_2-M_L\bar v_2\\
\dot{\bar v}_2&=-\alpha M_L {\bar v}_2+\alpha^2 M_L \bar r_2-\alpha^2 M_2^\top ({\bm \Pi}+{\bm \Pi}_1) 
\end{align*}
with  $\dot{\bar v}_1=0$ and ${\bm \Pi}_1\triangleq \nabla \tilde f(y)-\nabla \tilde f(r)$. By Assumption \ref{ass:convexity-strong}, ${\bm \Pi}$ and ${\bm \Pi}_1$ are   $\bar l$-Lipschitz with respect to $\bar r$ and $y-r$, respectively. By definitions,  $\dot{r}=M_1\dot{\bar r}_1+M_2\dot{\bar r}_2$.  Thus, there exist two constants $\bar c_1,\,\bar c_2>0$ such that  $||\dot{r}||^2< \bar c_1 W_{\rm o}(\bar r,\,\bar v_2)+\bar c_2||y-r||^2$. 

Using similar arguments as in the proof of Theorem \ref{thm:main},  we take the time derivative of $V_i$ and obtain
\begin{align*}
\dot{V}_i%&\leq 2 z_i^\top P_i [\frac{1}{\epsilon}A_{i1}z_i+ \frac{1}{\epsilon}A_{i2}\zeta_i+E_{i1}\dot{r}_i]\nonumber \\ 
%&\quad + \zeta_i(\frac{1}{\epsilon}A_{i3} z_i+\frac{1}{\epsilon}A_{i4} \zeta_i- \epsilon^{n_i-1}b_i\bar{\mathcal{N}}(\theta_i)\zeta_i+E_{i2}\dot{r}_i)\\
&\leq -\frac{2}{\epsilon}||z_i||^2+ \frac{1}{3\e}||z_i||^2+\frac{3}{\e}||P_i A_{i2}||^2\zeta_i^2+\frac{1}{3}||z_i||^2 \nonumber\\ 
&\quad +{3}||P_iE_{i1}||^2\dot{r}_i^2+\frac{1}{3\e}||z_i||^2+\frac{3}{\e}||A_{i3}||^2\zeta_i^2+\frac{1}{\e} A_{i4} \zeta_i^2 \nonumber\\ 
&\quad+ \epsilon^{n_i-1} b_i\bar{\mathcal{N}}(\theta_i)\zeta_i^2+ \zeta_i^2+||E_{i2}||^2\dot{r}_i^2 \nonumber \\ 
&=-(\frac{4}{3\e}-\frac{1}{3})||z_i||^2+( \epsilon^{n_i-1} b_i\bar{\mathcal{N}}(\theta_i)+\bar C_{i \theta_1})\zeta_i^2+ \bar C_{i\theta_2}\dot{r}_i^2 %\nonumber\\
%&=-\frac{1}{\e}||z_i||^2-( \epsilon^{n_i-1}  b_i\bar{\mathcal{N}}(\theta_i)-C_{i \theta_1})\dot{\theta}_i+ C_{i\theta_2}\dot{r}_i^2  %&\leq -(b_i\bar{\mathcal{N}}(\theta_i)-C_{i\theta_1})\dot{\theta}_i+ C_{i \theta_2}\dot{r}_i^2
\end{align*}
with $\bar C_{i\theta_1}=\frac{1}{\e}(3||P_i A_{i2}||^2+3||A_{i3}||^2+A_{i4})+1$ and $\bar C_{i \theta_2}=3||P_i E_{i1}||^2+||E_{i2}||^2$. Different from the proof of Theorem \ref{thm:main}, we avoid $\e$ here in handling the cross terms with $\dot{r}_i$ in order to dominate them %terms including $||{\dot r}||^2$ 
by deceasing $\e$. 

Denote $z=\mbox{col}(z_1,\,\dots,\,z_N)$ and $\hat {\bar C}_{\theta_2}=\max_{i\in \mathcal{N}}\{\bar C_{i\theta_2}\}$. We let ${\bar V}=\sum_{i=1}^N V_i+\sigma W_{\rm o}$ with $\sigma>0$ to be specified later. Its time derivative along the trajectory of the error system satisfies
\begin{align*}
\dot{\bar V}&\leq  \sum_{i=1}^N[-(\frac{4}{3\e}-\frac{1}{3})||z_i||^2+( \epsilon^{n_i-1} b_i\bar{\mathcal{N}}(\theta_i)+\bar C_{i \theta_1})\zeta_i^2+ \bar C_{i\theta_2}\dot{r}_i^2]\\
&\quad+\sigma \{-\frac{1}{2}W_{\rm o}-\alpha\bar r^\top [M_1~M_2]^\top {\bm \Pi}_1-\alpha^2  \bar v_2^\top M_2^\top {\bm \Pi}_1\}\\
%&\leq -(\frac{4}{3\e}-\frac{1}{3})||z||^2+\sum_{i=1}^N ( \epsilon^{n_i-1} b_i\bar{\mathcal{N}}(\theta_i)+\bar C_{i \theta_1})\zeta_i^2+ \hat {\bar C}_{\theta_2} ||\dot{r}||^2\\
%&\quad + \sigma[-\frac{1}{2}W_{\rm o}+\frac{1}{8}||\bar r||^2+\frac{1}{8\alpha^3}||\bar v_2||^2+(2\alpha^2+2\alpha^7) ||{\bm \Pi}_1||^2]\\
%&\leq -(\frac{4}{3\e}-\frac{1}{3})||z||^2+ \sum_{i=1}^N ( \epsilon^{n_i-1} b_i\bar{\mathcal{N}}(\theta_i)+\bar C_{i \theta_1})\zeta_i^2+ \hat {\bar C}_{\theta_2}\bar c_1 W_{\rm o}\\
%&\quad +\hat {\bar C}_{\theta_2}\bar c_2 ||y-r||^2-\frac{\sigma}{4}W_{\rm o}+2\sigma\bar l (\alpha^2+\alpha^7)||y-r||^2\\
&\leq -[\frac{4}{3\e}-\frac{1}{3}-2\sigma\bar l (\alpha^2+\alpha^7)-\hat {\bar C}_{\theta_2}\bar c_2]||z||^2-(\frac{\sigma}{4}- \hat {\bar C}_{\theta_2}\bar c_1)W_{\rm o}\\
&\quad +\sum_{i=1}^N ( \epsilon^{n_i-1} b_i\bar{\mathcal{N}}(\theta_i)+\bar C_{i \theta_1})\zeta_i^2  
\end{align*}
Letting  $\e^{\star}=\frac{1}{2\sigma\bar l (\alpha^2+\alpha^7)+\hat {\bar C}_{\theta_2}\bar c_2+1}$,  $\sigma \geq 8\max\{1,\,\hat {\bar C}_{\theta_2}\bar c_1\}$, and $0<\e<\e^{\star}$  gives
\begin{align*}
\dot{\bar V}&\leq -||z||^2- W_{\rm o}+  \sum_{i=1}^N ( \epsilon^{n_i-1} b_i\bar{\mathcal{N}}(\theta_i)+\bar C_{i \theta_1})\dot{\theta}_i  
\end{align*}
Recalling Lemma 4.4 in \cite{chen2019nussbaum}, one concludes the boundedness of $\bar V(t)$ and $\theta_i(t)$ over $[0,\,\infty)$. Thus, we can confirm the boundedness of all trajectories. Moreover, $||z||^2$ and $W_{\rm o}$ is integrable over $[0,\,\infty)$. By Lemma 8.2 in \cite{khalil2002nonlinear}, we have $z(t)\to 0$ and $W_{\rm o}\to 0$ as $t$ goes to $\infty$.  The rest proof can be complete by the same arguments as in the proof of Theorem \ref{thm:main}.
\pe

\begin{rem}\label{rem:unknown}
	In contrast with most optimal consensus works, multiple Nussbaum gains are employed in our proposed controllers \eqref{ctr:main} and \eqref{ctr:main-online}  to overcome the technical difficulties brought by unknown control directions. The obtained results definitely extend existing optimal consensus conclusions in \cite{zhang2017distributed,qiu2019distributed,tang2019cyb} to allow such type of system uncertainties.
\end{rem}

\begin{rem}\label{rem:consensus}
	Compared with the previous consensus results for multi-agent systems with or without unknown control directions in \cite{ren2008distributed, peng2014cooperative,tang2015output, huang2018fully,rezaei2018adaptive}, an optimization requirement is further considered in our formulation. Moreover, by letting $f_i(y)=(y-y_i(0))^2$, these two theorems provide an alternative way to achieve an average consensus goal  even these agents have unknown control directions.
\end{rem}

\section{Simulation}\label{sec:simu}

In this section, we propose two numerical examples to verify the effectiveness of our previous designs. 

\begin{figure}
\centering
\scalebox{0.6}{
	\centering
	\begin{tikzpicture}[node distance=1.8 cm, shorten >=1pt,  >=stealth',
	every state/.style ={circle, radius=5mm}] %minimum width=0.4cm}]
	\node[align=center,state](node1) {1};
	\node[align=center,state](node2)[right= of node1]{2};
	\node[align=center,state](node3)[right= of node2]{3};
	\node[align=center,state](node4)[above right= 0.3 and 1.2 of node3]{4};
	\node[align=center,state](node5)[above left = 0.3 and 1.2  of node4]{5};
	\node[align=center,state](node6)[left= of node5]{6};
    \node[align=center,state](node7)[left= of node6]{7};
	\node[align=center,state](node8)[above left= 0.3 and 1.2 of node1]{8};
	\path  (node1) edge [<-]  (node2)
	(node2) edge [<-]  (node3)
	(node3) edge [->]  (node4)
	(node4) edge [->]  (node5)
	(node5) edge [->]  (node6)
	(node6) edge [->]  (node7)
	(node7) edge [->, bend  left=20]  (node8)
	(node7) edge [<-,bend right=20]  (node8)
	(node8) edge [->, bend left=20]  (node1)
	(node8) edge [<-,bend right=20]  (node1)
	(node3) edge [->]  (node4)
	%(node1) edge[<-]  (node7)
	(node3) edge [<-]  (node5)
	(node1) edge [->,bend left=25]  (node3)
	(node5) edge [<-,bend left=25]  (node7)
	;
	\end{tikzpicture}
}
	\caption{Interconnection graph $\mathcal G$ in our examples.}\label{fig:graph}
\end{figure}
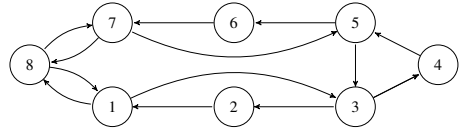

{\em Example 1}.  Consider an eight-agent network and each agent is described by double-integrator dynamics, that is,
\begin{align*} 
	\ddot{y}_i =b_i u_i,\quad i=1,\,\dots,\, 8
\end{align*}
Assume their interconnection topology is depicted in Fig.\ref{fig:graph} with unity weights. Assumption \ref{ass:graph} can be verified.  We are going to solve an average consensus for these agents. 

\begin{figure}
	\centering
	\includegraphics[width=0.38\textwidth]{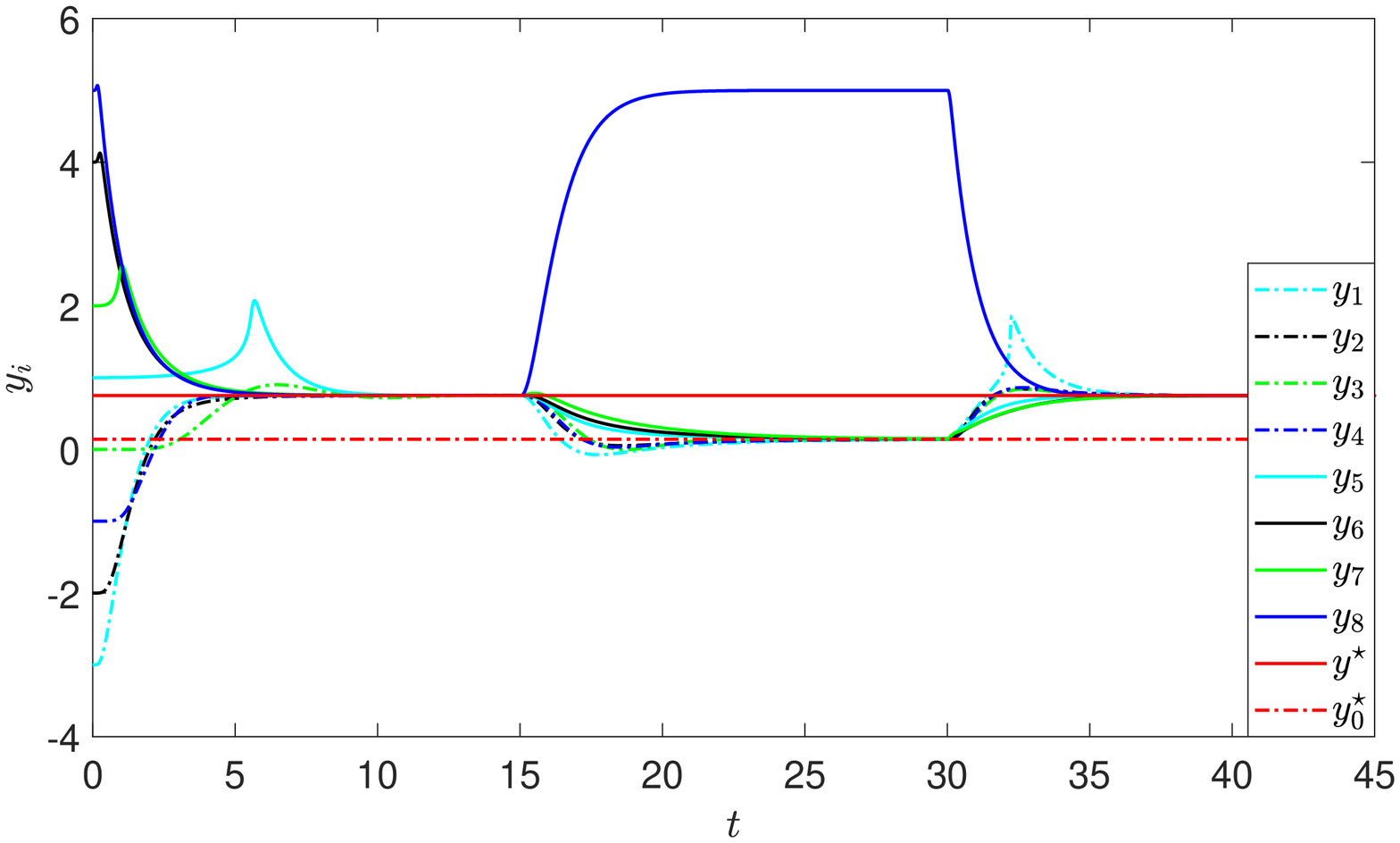}\\ %%% 0.80
	\caption{Profiles of agent output $y_i(t)$ in {\em Example 1}.}\label{fig:average-ouptut}
\end{figure}

According to Remark \ref{rem:consensus}, we let $f_i(y)=(y-y_i(0))^2$ for $i=1,\,\dots,\,8$ and use the controller \eqref{ctr:main} with $n_i=2$ to complete the design. For simulation, we set $b_1=\dots=b_4=-1$, $b_5=\dots=b_8=1$, and $y(0)=[-3~-2~0~-1~1~4~2~5]^\top$. Distributed controller \eqref{ctr:main} with $\e=1$, $k_{i1}=1$ for $i=1,\,\dots,\,8$, and $\bar{\mathcal{N}}(\theta)=\theta^2\sin \theta$ is then applied to solve this problem.  To make it more interesting, we cut all links associated with node $8$ at $t=15{\rm s}$ and then add them back at $t=30{\rm s}$. %The profiles of agents' control efforts and adaptive gains are showed in Figs.~\ref{fig:average-control} and \ref{fig:average-gain}, respectively. 
The simulation result is depicted in Fig.~\ref{fig:average-ouptut}. At first, the outputs of agents are observed to reach an average consensus on  $y^{\star}=\frac{\sum_{i=1}^8y_i(0)}{8}=0.75$. Then, $y_8(t)$ converges to its local optimizer $y_8(0)=5$ while the other agents reach a consensus on $y^{\star}_0=\frac{\sum_{i=1}^7y_i(0)}{7}=0.143$. After the links are added back, the average consensus for all agents is quickly recovered at $y^{\star}$. This verifies the robustness of our algorithms enabling plug-and-play operations.

%\begin{figure}
%	\centering
%	\subfigure[]{
%		\includegraphics[width=0.22\textwidth]{opt-adaptive-output-opt.eps}
%	} 
%	\subfigure[]{
%		\includegraphics[width=0.22\textwidth]{opt-adaptive-gain-average.eps}
%	}
%	\caption{Profiles of agent output $y_i(t)$ and control effort $u_i(t)$ in {\em Example 2} under the controller \eqref{ctr:offline}.}\label{fig:opt}
%\end{figure}

{\em Example 2}.  Consider the optimal consensus problem for a heterogeneous multi-agent system with agents described by $${y}_i^{(m_i)}=b_iu_i,\quad i=1,\,\dots,\,8$$ with the same topology as that in {\em Example 1}. Here, $m_1=m_5=1$, $m_2=m_6=2$, $m_3=m_7=3$, and $m_4=m_8=4$.

The local cost functions are taken as ${f_1}(y) =f_5(y)= (y-8)^2$, ${f_2}(y) =f_6(y) =\frac{y^2}{{20\sqrt {y^2 + 1} }} + y^2$, ${f_3}(y) = f_7(y) =\frac{y^2}{80\ln {({y^2} + 2} )} + (y -5)^2$, ${f_4}(y) =f_8(y)=  \ln \left( {{e^{ - 0.05{y}}} + {e^{0.05{y}}}} \right) + y^2$. 
Assumption \ref{ass:convexity-strong} is confirmed with $\underline{l}=1$, $\bar{l}=3$ as that in \cite{tang2020optimal}. Moreover, the global optimal point can be obtained numerically as $y^{\star}=3.24$. Since these agents are of heterogeneous orders and unknown high-frequency gains, the rules developed in \cite{zhang2017distributed,tang2019cyb} fail to tackle this problem. Nevertheless, according to Theorems \ref{thm:main} and \ref{thm:main-online}, we can utilize controller \eqref{ctr:main} or \eqref{ctr:main-online} to solve it. 

\begin{figure}
	\centering
	\includegraphics[width=0.37\textwidth]{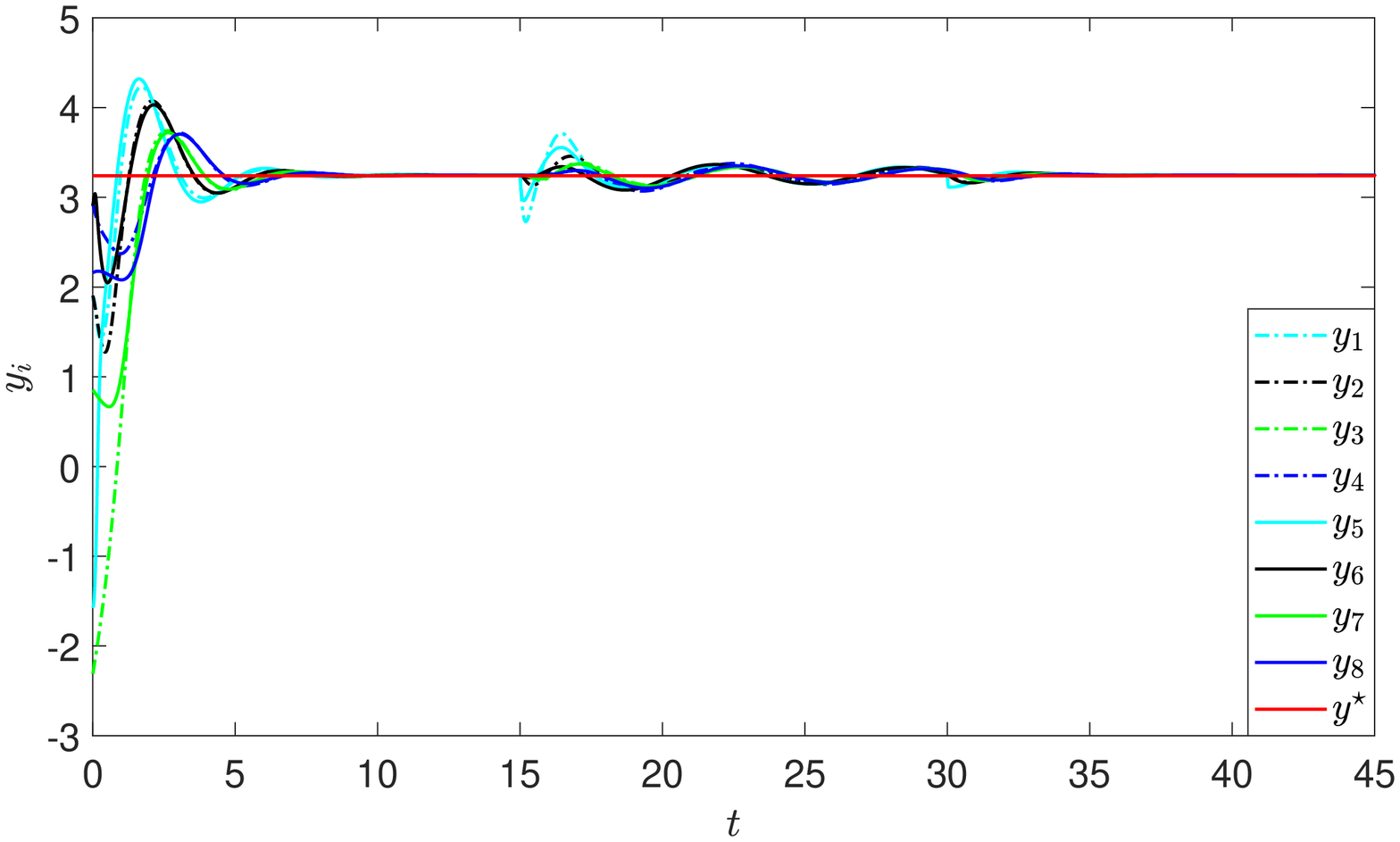}\\ %%% 0.80
	\caption{Profiles of control effort $y_i(t)$ in {\em Example 2}.}\label{fig:opt-ouptut}
\end{figure}

\begin{figure}
	\centering
	\includegraphics[width=0.44\textwidth]{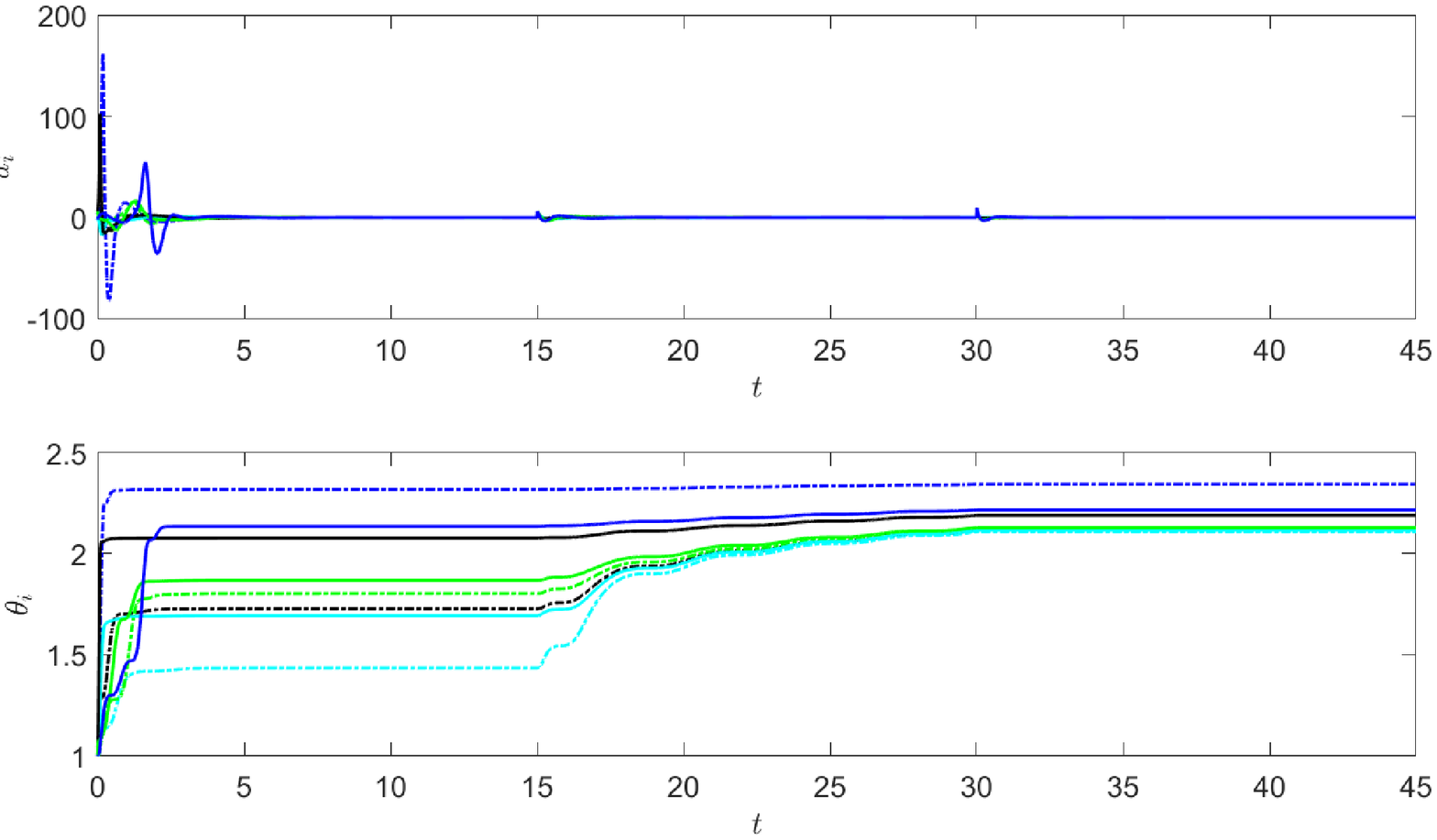}\\ %%% 0.80
	\caption{Profiles of control efforts $u_i(t)$ and adaptive gain $\theta_i(t)$ in {\em Example 2}.}\label{fig:opt-control}
\end{figure}

For simulation, we let $b_1=\dots=b_8=-1$. Choose $k_{21}=k_{61}=1$, $k_{31}=k_{71}=1$, $k_{32}=k_{72}=2$, $k_{41}=k_{81}=1$, $k_{42}=k_{82}=3$, $k_{43}=k_{83}=3$, $\epsilon=0.5$, and $\bar{\mathcal{N}}(\theta)=e^{\theta^2}\sin \theta$ for controller \eqref{ctr:main-online}.  To verify the robustness of our algorithm, we add an actuated disturbance $10\sin(t)$ for all agents during $15{\rm s}\leq t\leq 30 {\rm s}$. %The profiles of agents' control efforts and adaptive gains are showed in Figs.~\ref{fig:average-control} and \ref{fig:average-gain}, respectively. 
The simulation result is depicted in Figs.~\ref{fig:opt-ouptut} and \ref{fig:opt-control}. One can observe that all agents quickly reach an optimal consensus on  $y^{\star}=3.24$ at first while the profiles of agents' control efforts are maintained bounded. Then, the expected exact optimal consensus is broken due to actuated disturbances but the error $|y_i-y^{\star}|$ is still bounded. These observations verify the efficacy and robustness of our adaptive optimal consensus algorithms in handling both heterogeneous agent dynamics and unknown control directions.

%\begin{figure}
%	\centering
%	\subfigure[]{
%		\includegraphics[width=0.22\textwidth]{opt-adaptive-output-opt.eps}
%		} 			
%	\subfigure[]{
%		\includegraphics[width=0.22\textwidth]{opt-adaptive-gain-average.eps}
%	}
%	\caption{Profiles of control efforts $u_i(t)$ and adaptive gain $\theta_i(t)$ in {\em Example 2}.}\label{fig:opt-control}
%\end{figure}

\section{Conclusion}\label{sec:con}

An optimal consensus problem has been discussed for a high-order multi-agent system without a prior knowledge of the control directions.  By an embedded design, we finally propose two Nussbaum-type distributed controllers to solve it under different information circumstances.   Further works will include improvement of transient performances and  extensions with more general agent dynamics.

%\section*{References}
%\bibliographystyle{elsarticle-num-names}
\bibliographystyle{IEEEtrans}
\bibliography{opt_nuss}
\end{document}